\documentclass{amsart}
\usepackage{verbatim}
\usepackage{amssymb}
\usepackage{amsbsy}
\usepackage{amscd}
\usepackage{amsmath}
\usepackage{amsthm}
\usepackage[mathscr]{eucal}

\newenvironment{NB}{
{\bf NB}. \footnotesize
}{}
\renewenvironment{NB}{
\comment
  }{\endcomment}

\theoremstyle{plain}
 \newtheorem{thm}{Theorem}[section]
 \newtheorem{lem}[thm]{Lemma}
 \newtheorem{prop}[thm]{Proposition}
 \newtheorem{cor}[thm]{Corollary}
 
\theoremstyle{definition}
 \newtheorem{defn}{Definition}[section]
\theoremstyle{remark}
 \newtheorem{rem}{Remark}[section]

\def\Bbb{\mathbb}
\def\frak{\mathfrak}
\def\cal{\mathcal}

\newcommand{\Ext}{\operatorname{Ext}}
\newcommand{\Hom}{\operatorname{Hom}}

\newcommand{\im}{\operatorname{im}}
\newcommand{\Aut}{\operatorname{Aut}}
\newcommand{\rk}{\operatorname{rk}}

\newcommand{\NS}{\operatorname{NS}}
\newcommand{\coker}{\operatorname{coker}}
\newcommand{\Pic}{\operatorname{Pic}}

\newcommand{\Alb}{\operatorname{Alb}}
\newcommand{\Hilb}{\operatorname{Hilb}}
\newcommand{\Quot}{\operatorname{Quot}}

\newcommand{\Coh}{\operatorname{Coh}}

\newcommand{\IT}{\operatorname{IT}}
\newcommand{\WIT}{\operatorname{WIT}}

\newcommand{\Div}{\operatorname{Div}}

\font\b=cmr10 scaled \magstep5
\def\bigzerou{\smash{\lower1.7ex\hbox{\b 0}}}

\numberwithin{equation}{section}
\begin{document}

\title
{
Fourier-Mukai transform on abelian surfaces
}
\author{K\={o}ta Yoshioka
}
 \address{Department of Mathematics, Faculty of Science,
Kobe University,
Kobe, 657, Japan}
\email{yoshioka@math.kobe-u.ac.jp}
 
\maketitle

\section{Introduction}\label{sect:intro}
Let $(X,H)$ be a pair of an abelian surface $X$ and an ample
line bundle $H$ on $X$.
Let $\langle\;\;,\;\;\rangle$ be a bilinear form on
$H^{ev}(X,{\Bbb Z}):=\bigoplus_i H^{2i}(X,{\Bbb Z})$
defined by
\begin{equation}
\langle x,y \rangle:=\int_X (x_1 \cup y_1-x_0 \cup y_2-x_2 \cup y_0)
\end{equation}
where $x=(x_0,x_1,x_2), y=(y_0,y_1,y_2)$ with 
$x_i,y_i \in H^{2i}(X,{\Bbb Z})$.
For an object $E \in {\bf D}(X)$,
we define the Mukai vector $v(E) \in H^{ev}(X,{\Bbb Z})$
of $E$ as the Chern character of $E$.
We also call an element $v \in H^*(X,{\Bbb Z})$
Mukai vector, if $v=v(E)$ for an object $E \in {\bf D}(X)$.

We denote the coarse moduli space of $S$-equivalence classes of
semi-stable sheaves $E$ with
$v(E)=v$ by $\overline{M}_H(v)$
and the open subscheme consisting of stable sheaves by $M_H(v)$.
We also denote the moduli stack of semi-stable sheaves by
${\cal M}_H(v)^{ss}$. 
Let $Y:=M_H(v_0)$ be the moduli space of stable semi-homogeneous sheaves
on $X$. Assume that $Y$ is a fine moduli space, that is,
there is a universal family ${\bf E}$ 
on $Y \times X$.
We define the integral functor $\Phi_{Y \to X}^{{\bf E}}$ as
\begin{equation}
\begin{matrix}
{\bf D}(Y) & \to & {\bf D}(X)\\
y & \mapsto &  {\bf R}p_{X*}({\bf E} \otimes p_{Y}^*(y)),
\end{matrix}
\end{equation}
where $p_{X}:X \times Y \to X$ (resp. $p_Y:X \times Y \to Y$)
be the projection.
Let ${\bf D}(X)_{op}$ be the opposite category of 
${\bf D}(X)$ and define the equivalence
\begin{equation}
\begin{matrix}
D:& {\bf D}(X) & \to & {\bf D}(X)_{op}\\
& x & \mapsto &  x^{\vee}={\bf R}{\cal H}om(x,{\cal O}_X).
\end{matrix}
\end{equation}

\begin{defn}
We call equivalences
${\bf D}(Y) \to {\bf D}(X)$
and ${\bf D}(Y) \to {\bf D}(X)_{op}$
the {\it Fourier-Mukai transform}.
\end{defn}
$\Psi_{Y \to X}^{{\bf E}}:H^*(Y,{\Bbb Z}) \to H^*(X,{\Bbb Z})$
denotes the cohomological transform induced by ${\bf E}$.

\begin{thm}\label{thm:main}
Let $w \in H^*(Y,{\Bbb Z})$ 
be a primitive Mukai vector with $\langle w^2 \rangle>0$. 
Let $H'$ be an ample divisor on $Y$ witch is general with respect to
$w$.
We set $v=\Psi_{Y \to X}^{{\bf E}}(w)$.
We assume that $H$ is general with respect to $v$.
Then there is an autoequivalence 
$\Phi_{X \to X}^{{\bf F}}:{\bf D}(X) \to {\bf D}(X)$ such that
for a general $E \in M_{H'}(v)$,
$F:=\Phi_{X \to X}^{{\bf F}} \circ \Phi_{Y \to X}^{{\bf E}}(E)$
is a stable sheaf with $v(F)=v$
or $F^{\vee}$ is a stable sheaf with $v(F^{\vee})=v$.
In particular, 
$M_{H'}(w)$ is birationally equivalent to 
$M_H(v)$.
\end{thm}
Since the moduli space of semi-stable sheaves is irreducible,
the same assertion in Theorem \ref{thm:main} 
also holds for a general stable sheaf
$E$ with a non-primitive vector.
This is a partial generalization of a result in \cite{Y:12},
which is conjectured in \cite[Conj. 4.16]{Y:7}.
If $X$ is a $K3$ surface, then 
a similar result is conjectured by Tom Bridgeland.
In particular, the idea of replacing 
$\Phi_{Y \to X}^{{\bf E}}(E)$ by 
another Fourier-Mukai transform
$\Phi_{X \to X}^{{\bf F}} \circ \Phi_{Y \to X}^{{\bf E}}(E)$ 
is due to him.

\section{Preliminaries}\label{sect:pre}

\subsection{A family of 2-extensions}
In this section, we recall or prepare some necessary results 
to prove Theorem \ref{thm:main}.
We start with a possibly well-known result
on a family of 2-extensions.
   
\begin{defn}
Let 
\begin{equation}
{\cal V}_{\bullet}:\cdots \to {\cal V}_{-1} \to {\cal V}_0 \to
{\cal V}_1 \to \cdots
\end{equation}
be a complex on $X \times T$.
${\cal V}_{\bullet}$ is flat, if ${\cal V}_i$ are flat over $T$.
\end{defn}

We shall construct a family of 2-extensions:
\begin{equation}
0 \to A_0 \to V_0 \to V_1 \to A_1 \to 0.
\end{equation}
Let $v_0,v_1$ be Mukai vectors of coherent sheaves on $X$.
Let $Q_i$, $i=0,1$ be the open subscheme of the quot-scheme
$\Quot_{W_i \otimes {\cal O}_X(-n_i)/X}^{v_i}$
parametrizing all quotients
$W_i \otimes {\cal O}_X(-n_i) \to A_i$ with $v(A_i)=v_i$ such that
$W_i \cong H^0(X,A_i(n_i))$ and 
$H^j(X,A_i(n_i))=0, j>0$.
Let ${\cal A}_i$ be the universal quotient and
${\cal I}_i$ the universal subsheaf.
We take an integer $m $ such that
(i) $R^j p_{Q_1 *}({\cal I}_1(m))=0, j>0$, 
(ii) ${\cal U}:=p_{Q_1 *}({\cal K}_1(m))$ is locally free
and (iii)
$\psi_0:p_{Q_1}^* ({\cal U}) \to {\cal I}_1(m)$
is surjective. 
We set
${\cal J}:=\ker(\psi_0)(-m)$. 
Let $\widetilde{Q}_1 \to Q_1$ be the principal $GL$-bundle
associated to ${\cal U}$.
Then we have a trivialization
${\cal U} \cong U \otimes {\cal O}_{\widetilde{Q}_1}$, where 
$U$ is a vector space.
Let $\widetilde{\cal I}_i$ (resp. $\widetilde{\cal J}, \widetilde{\cal A}_i$)
be the pull-back of ${\cal I}_i$ (resp. ${\cal J},{\cal A}_i$)
to $Q_0 \times \widetilde{Q}_1 \times X$.
Then we have exact sequences
\begin{equation}
0 \to \widetilde{\cal J} \to 
U \otimes {\cal O}_{Q_0 \times \widetilde{Q}_1 \times X}(-m)
\to \widetilde{\cal I}_1 \to 0,
\end{equation}

\begin{equation}
0 \to \widetilde{\cal I}_i \to W_i \otimes 
{\cal O}_{Q_0 \times \widetilde{Q}_1 \times X}(-n_i)
\to \widetilde{\cal A}_i \to 0.
\end{equation}
If $m$ is sufficiently large, then 
$\Ext^j_{p_{Q_0 \times Q_1}}(\widetilde{\cal J},\widetilde{\cal A}_0)=0$
and
${\Bbb E}:=
\Hom_{p_{Q_0 \times Q_1}}(\widetilde{\cal J},\widetilde{\cal A}_0)$ 
is locally free.
We have an exact sequence:
\begin{equation}
0 \to \Hom_{p_{Q_0 \times \widetilde{Q}_1}}
(\widetilde{\cal I}_1,\widetilde{\cal A}_0) \to
\Hom_{p_{Q_0 \times \widetilde{Q}_1}}
(U \otimes {\cal O}_{Q_0 \times \widetilde{Q}_1 \times X}(-m),
\widetilde{\cal A}_0) \to {\Bbb E} \to
\Ext^1_{p_{Q_0 \times \widetilde{Q}_1}}
(\widetilde{\cal I}_1,\widetilde{\cal A}_0) \to 0.
\end{equation}
Let $\pi:P \to Q_0 \times \widetilde{Q}_1$ 
be the associated vector bundle
of ${\Bbb E} $. Then we have a family of extensions
\begin{equation}
0 \to (\pi \times 1_X)^*(\widetilde{\cal A}_0) \to
\widehat{\cal V}_0 \to
(\pi \times 1_X)^*(\widetilde{\cal I}_1) \to 0.
\end{equation}
We set $\widehat{\cal V}_1:=W_1 \otimes {\cal O}_{P \times X}(-n_1)$.
Then we have a family of 
complexes
\begin{equation}
\widehat{\cal V}_{\bullet}:\widehat{\cal V}_0 \to \widehat{\cal V}_1
\end{equation}
such that
$\widehat{\cal V}_i$ are flat over $P$,
$H^j(\widehat{\cal V}_{\bullet})$ are flat over $P$
and $H^j(\widehat{\cal V}_{\bullet})_{x}=({\cal A}_j)_{\pi(x)}$.

Let $S_i$ be a bounded set of coherent sheaves $E_i$ on $X$
with the Mukai vector $v_i$.

\begin{prop}
Let ${\cal V}_{\bullet}$ be a $T$-flat family of complexes on $X$
parametrized by $T$
such that $H^i({\cal V}_{\bullet})$ are flat families
of coherent sheaves belonging to $S_i$.
Then for any point $t \in T$,
there is a neighborhood $T_0$ of $t$ with the following prpperty:
there is a quasi-isomorphism 
${\cal V}_{\bullet}' \to {\cal V}_{\bullet}$
and a morphism $f:T \to P$ such that $f^*(\widehat{\cal V}_{\bullet})$ is 
quasi-isomorphic to ${\cal V}_{\bullet}'$.
\end{prop}

\begin{proof}
Construction of ${\cal V}_{\bullet}' \to {\cal V}_{\bullet}$:
Let ${\cal V}_{\bullet}:=({\cal V}_0 \overset{\phi}{\to} {\cal V}_1)$ 
be a flat family of 
complexes on $X \times T$ 
such that $H^i({\cal V}_{\bullet})$ are flat over $T$.
Let ${\cal B}$ be the kernel of 
${\cal V}_1 \to H^1({\cal V}_{\bullet})$.
For a sufficiently large $n$, 
$R^j p_{T *}({\cal B}(n))=R^j p_{T *}({\cal V}_1(n))=
R^j p_{T *}(H^1({\cal V}_{\bullet})(n))=0$ for $j>0$ and
we have an exact and commutative diagram:
\begin{equation}
\begin{CD}
0 @>>> p_{T }^*(p_{T *}({\cal B}(n))) @>>>
p_{T }^*(p_{T *}({\cal V}_1(n))) @>>>
p_{T }^*(p_{T *}(H^1({\cal V}_{\bullet})(n))) @>>>0\\
@. @VVV @VVV @VV{\psi}V\\
0@>>> {\cal B}(n)@>>>
{\cal V}_1(n)@>>>
H^1({\cal V}_{\bullet})(n)@>>> 0\\
@. @VVV @VVV @VVV\\
@. 0 @. 0@. 0
\end{CD}
\end{equation} 
By shrinking $T$ if necessary,, we may assume that
there is a lifting 
$\widetilde{\psi}:p_{T }^*(p_{T *}(H^1({\cal V}_{\bullet})(n))) \to
{\cal V}_1(n)$ of 
$\psi$.
We set ${\cal V}_1':=p_{T }^*(p_{T *}(H^1({\cal V}_{\bullet})(n)))(-n)$
and set ${\cal K}_1:=\ker (\psi)(-n)$.
Then we have a homomorphism ${\cal K}_1 \to {\cal B}$.
Let ${\cal V}_0'$ be a coherent sheaf fitting in the diagram

\begin{equation}
\begin{CD}
0 @>>>H^0({\cal V}_{\bullet})  @>>> {\cal V}_0' @>>> {\cal K}_1 @>>> 0\\
@. @| @VVV @VVV @.\\
0 @>>> H^0({\cal V}_{\bullet}) @>>> {\cal V}_0 @>>> {\cal B} @>>> 0.
\end{CD}
\end{equation}
Then ${\cal V}_{\bullet}':{\cal V}_0' \to {\cal V}_1'$ is quasi-isomorphic
to ${\cal V}_{\bullet}$.
We shall show that there is a local morphism $f:T \to P$
with a quasi-isomorphism
${\cal V}_{\bullet}' \to (f \times 1_X)^*(\widehat{\cal V}_{\bullet})$
for a sufficiently large $n$.

Construction of $f:T \to P$:
We take a trivialization $p_{T *}(H^i({\cal V}_{\bullet})(n_i)) \cong
W_i \otimes {\cal O}_T$.
Then we have a morphism $h:T \to Q_0 \times Q_1$
such that $(h \times 1_X)^*({\cal A}_i)
=H^i({\cal V}_{\bullet})$ as quotients of
$W_i \otimes {\cal O}_{T \times X}(-n_i)$.
If $n$ is sufficiently large, then
$\Hom({\cal V}_1',W_1 \otimes {\cal O}_{T \times X}(-n_1))
\to \Hom({\cal V}_1',H^1({\cal V}_{\bullet}))$ is surjective.
Hence there is a homomorphism 
${\cal V}_1' \to W_1 \otimes {\cal O}_{T \times X}(-n_1)$
and a commutative diagram
\begin{equation}
\begin{CD}
0 @>>> {\cal K}_1 @>>> {\cal V}_1' @>>> H^1({\cal V}_{\bullet}) @>>> 0\\
@. @VVV @VVV @| @. \\
0 @>>> (h \times 1_X)^*({\cal I}_1) @>>> 
W_1 \otimes {\cal O}_{T \times X}(-n_1) 
@>>> H^1({\cal V}_{\bullet})@>>> 0,
\end{CD}
\end{equation}
where ${\cal I}_1$ means the pull-back of ${\cal I}_1$ to
$Q_0 \times Q_1 \times X$.
By our choice of $n_i$ and $n$, we have
\begin{equation}
\begin{CD}
\Ext^1_{p_{T}}({\cal K}_1,H^0({\cal V}_{\bullet})) & \;\cong \;& 
\Ext^2_{p_{T}}(H^1({\cal V}_{\bullet}),H^0({\cal V}_{\bullet}))\\
@AAA @| \\
\Ext^1_{p_{T}}((h \times 1_X)^*({\cal I}_1),H^0({\cal V}_{\bullet})) 
& \;\cong \;& 
\Ext^2_{p_{T}}(H^1({\cal V}_{\bullet}),H^0({\cal V}_{\bullet})).
\end{CD}
\end{equation}
Shrinking $T$ if necessary,
there is a coherent sheaf $\widetilde{\cal V}_0$ on $T \times X$
fitting in 
the following diagram:
\begin{equation}
\begin{CD}
0 @>>>H^0({\cal V}_{\bullet})  @>>> {\cal V}_0' @>>> {\cal K}_1 @>>> 0\\
@. @| @VVV @VVV @.\\
0 @>>> H^0({\cal V}_{\bullet}) @>>> \widetilde{\cal V}_0 @>>> 
(h \times 1_X)^*({\cal I}_1) @>>> 0.
\end{CD}
\end{equation}
Then by shrinking $T$, we have a morphism
$f:T \to P^s$ with a commutative diagram:
\begin{equation}
\begin{CD}
0 @>>>H^0({\cal V}_{\bullet})  @>>> (f \times 1_X)^*(\widehat{\cal V}_0) 
@>>> (f \times 1_X)^*(\widetilde{\cal I}_1) @>>> 0\\
@. @| @VVV @| @.\\
0 @>>> H^0({\cal V}_{\bullet}) @>>> \widetilde{\cal V}_0 @>>> 
(f \times 1_X)^*(\widetilde{\cal I}_1) @>>> 0.
\end{CD}
\end{equation}
Therefore $(f \times 1_X)^*(\widehat{\cal V}_{\bullet})
\cong (\widetilde{\cal V}_0 \to W_1 \otimes {\cal O}_{T \times X}(-n_1))$
is quasi-isomorphic to ${\cal V}_{\bullet}'$.
\end{proof}

\subsection{Albanese map}\label{subsect:albanese}
Let $\widehat{X}$ be the dual abelian variety of $X$ and
${\bf P}$ the Poincar\'{e} line bundle on $\widehat{X} \times X$.
Let ${\frak a}:{\bf D}(X) \to \Pic(\widehat{X}) \times \Pic(X)$
be the morphism sending $E$ to
$(\det \Phi_{X \to \widehat{X}}^{{\bf P}}(E),\det(E))$.
For a family of coherent sheaves ${\cal E}$
parametrized by a connected scheme $T$, 
we also have a morphism ${\frak a}:T \to
X \times \widehat{X}$ (up to translation).  

We quote the following assertions from 
\cite[Thm. 0.1, Lemma 4.3, Prop. 4.4]{Y:7}.
\begin{prop}\label{prop:albanese}
Let $v$ be a Mukai vector.
\begin{enumerate}
\item
Let ${\cal E}$ be a family of coherent sheaves on $X$ with
$v({\cal E}_q)=v$ parametrized by a scheme $Q$.
Assume that for any point $(x,y) \in X \times \widehat{X}$,
$T^*_x({\cal E}_q) \otimes {\bf P}_y \cong 
{\cal E}_{q'}$ for a point $q' \in Q$.
Then $\dim {\frak a}(Q) \geq 2$ and $\dim {\frak a}(Q)=4$ if
$\langle v^2 \rangle>0$.
\item
If $v$ is a primitive Mukai vector with $\langle v^2 \rangle>0$.
Then ${\frak a}:M_H(v) \to \widehat{X} \times X$ is the Albanese map.
\end{enumerate}
\end{prop}
In the case where $\langle v(E)^2 \rangle=0$, 
we use Lemma 4.3 in \cite{Y:7} and
the fact that $\phi_L=0$ if and only if $c_1(L)=0$.

%

\section{Proof of Theorem \ref{thm:main}}

\subsection{Fourier-Mukai transform of a general stable sheaf}
\label{subsect:FM}
Let $Y$ be a moduli space of stable semi-homogeneous sheaves on $X$.
Assume that there is a universal family ${\bf E}$ on $Y \times X$.
Then we have a Fourier-Mukai transform
$\Phi_{Y \to X}^{{\bf E}}:{\bf D}(Y) \to {\bf D}(X)$.
If $\dim {\bf E}_y=0$, $y \in Y$, then
$\Phi_{Y \to X}^{{\bf E}}$ comes from an equivalence
$\Coh(Y) \to \Coh(X)$.
This case is easier to treat than other cases.
In particular, the proof of Theorem \ref{thm:main}
is reduced to the case treated in \ref{subsect:sheaf}.
Hence we assume that $\dim {\bf E}_y \geq 1$, $y \in Y$.

\begin{thm}\label{thm:complex}
Let $w$ be a primitive Mukai vector such that
$\langle w^2 \rangle>0$. 
If $\Phi_{Y \to X}^{{\bf E}}(E)$ is not a sheaf up to shift
for all $E \in M_{H'}(w)$, then
there is an integer $k$ such that
for a general $E \in M_{H'}(w)$,
$\Phi_{Y \to X}^{{\bf E}}(E)[k]$ fits in an exact triangle
\begin{equation}
A_0 \to \Phi_{Y \to X}^{{\bf E}}(E)[k] \to A_1[-1] \to A_0[1],
\end{equation}
where $A_i, i=0,1$ are semi-homogeneous sheaves of $v(A_i)=n_i' v_i'$,
$(n_0'-1)(n_1'-1)=0$ and $\langle v_0',v_1' \rangle=-1$.
In particular
$\Phi_{{\bf E}}$ induces a birational map
$M_{H'}(w) \cdots \to M_H(v)$,
if $\NS(X) \cong {\Bbb Z}$ and
$v \ne \pm (v_0'-n v_1')$ for all isotropic vectors $v_0', v_1'$ with
$\langle v_0',v_1' \rangle=-1$,
where $v:=\Psi_{Y \to X}^{{\bf E}}(w)$.
\end{thm}

\begin{proof}
Let $Q(w)$ be the open subscheme of $\Quot_{{\cal O}_Y(-n)^{\oplus N}/Y}^w$
such that $M_{H'}(w)$ is the geometric quotient of
$Q(w)$ by the action of $PGL(N)$.
Let ${\cal E}$ be the universal family on
$Q(w) \times Y$.
Then for a point $q \in Q(w)$, we have
\begin{equation}
\begin{cases}
H^0(\Phi_{Y \to X}^{{\bf E}}({\cal E}_q))=0, & 
\mu({\cal E}_q \otimes {\bf E}_x) \leq 0\\
H^2(\Phi_{Y \to X}^{{\bf E}}({\cal E}_q))=0, & 
\mu({\cal E}_q \otimes {\bf E}_x) > 0,
\end{cases}
\end{equation}
where $x \in X$.
Hence $\Phi_{Y \to X}^{{\bf E}}({\cal E})[k]$ is a family of complexes
represented by
\begin{equation}
{\cal V}_{\bullet}:{\cal V}_0 \to {\cal V}_1,
\end{equation}
where $k=1$ or $k=0$.
Assume that $\WIT$ does not hold for all ${\cal E}_q$.
Let $Q(w)_0$ be the open subscheme such that
$H^i({\cal V}_{\bullet})$ are flat over $Q(w)_0$.
Let $S_i$ be the bounded set of coherent sheaves
$H^i({\cal V}_{\bullet})_q, q \in Q(w)_0$.
We set $v_i:=v(H^i({\cal V}_{\bullet})_q)$ and consider the family of
complexes $\widehat{\cal V}_{\bullet}$ parametrized by $P$.
Then for any point $q \in Q(w)_0$, there is a neighborhood
$Q_q$ of $q$ and a morphism $f_q:Q_q \to P$.
We note that $Q(w)_0$ is $GL(N)$-invariant.
We set $M_H(w)_0:=Q(w)_0/GL(N)$. 
By shrinking $Q(w)_0$ to an open subscheme, 
we may assume that
the Harder-Narasimhan filtrations
$0 \subset F_i^1 \subset F_i^2 \subset \cdots \subset 
F_i^{s_i}=H^i({\cal V}_{\bullet})_q, q \in Q(w)_0$ 
form a flat family of filtrations over $Q(w)_0$, that is,
$F_i^j/F_i^{j-1}$ form flat families of sheaves.
We set $v_i^j:=v(F_i^j/F_i^{j-1})$ and 
consider the locally closed subset $Q_i'$
of $Q_i$ such that
\begin{equation}
Q_i'=\left\{
W_i \otimes {\cal O}_X(-n_i) \to A_i \left|
\begin{split}
& \text{ the Harder-Narasimhan filtration of $A_i$ is }\\
& \text{ $0 \subset F_i^1 \subset F_i^2 \subset \cdots \subset 
F_i^{s_i}=A_i$, $v(F_i^j/F_i^{j-1})=v_i^j$}
\end{split}
\right.
\right\}
\end{equation}
(cf. Remark \ref{rem:component}).
Then we have a morphism 
\begin{equation}
\begin{matrix}
{\frak a}_i':&
Q_i' & \to & \prod_j \overline{M}_H(v_i^j) & 
\to & (X \times \widehat{X})^{s_i}\\
& A_i & \mapsto & \prod_j F_i^j/F_i^{j-1}& \mapsto&
 \prod_j {\frak a}(F_i^j/F_i^{j-1}).
\end{matrix}
\end{equation} 
By the proof of \cite[Thm. 4.14]{Y:7}, 
we get $\dim {\frak a}_i'(Q_i') \geq 2s_i$.
Indeed if $n_i$ is sufficiently large, then
we can show that the quotient stack 
$[Q_i'/GL(W_i)]$ is an affine bundle
over $\prod_j {\cal M}_H(v_i^j)^{ss}$
(see \cite[sect.2.2, in particular Prop. 2.5]{Y:11}). 
Combining this with Proposition \ref{prop:albanese},
we get $\dim {\frak a}_i'(Q_i') \geq 2s_i$.
We set $P':=P \times_{(Q_0 \times Q_1)} (Q_0' \times Q_1')$.
Then the image of $f_q:Q_q \to P$ is contained in $P'$.
Let 
${\frak b}:P' \to Q_0' \times Q_1' \to (X \times \widehat{X})^{s_0+s_1}$
be the morphism defined by 
composing $\pi$ with ${\frak a}_0' \times {\frak a}_1'$.
Then $\dim \im {\frak b} \geq 4$ and if the equality holds,
then $s_0=s_1=1$ and $\langle v_0^2 \rangle=\langle v_1^2 \rangle=0$.
Thus $Q_i'$ are open subset of $Q_i$ and  
${\cal A}_i$ are families of semi-homogeneous sheaves.

Let $P^s$ be the open subset of $P'$
such that $\Phi_{X \to Y}^{{\bf E}^{\vee}}(\widehat{\cal V}_{\bullet})[2-k]$
is a family of stable sheaves.
Then we have a morphism $g:P^s \to M_H(w)$.
Obviously $g \circ f_q:Q_q \to M_H(w)$ is the restriction of the
quotient map $\varpi$.
Combining with ${\frak b}$, we have a morphism
$Q_q \to P^s \to (X \times \widehat{X})^{s_0+s_1}$.
This is the morphism determined by ${\cal E}_q$:
\begin{equation}
Q(w)_0 \ni q \mapsto 
({\frak a}_0'(H^0(\Phi_{Y \to X}^{{\bf E}}({\cal E}_q)[k])),
{\frak a}_1'(H^1(\Phi_{Y \to X}^{{\bf E}}({\cal E}_q)[k]))) \in
(X \times \widehat{X})^{s_0} \times (X \times \widehat{X})^{s_1}.
\end{equation}
Hence we have a morphism 
${\frak c}:M_H(v)_0 \to (X \times \widehat{X})^{s_0+s_1}$
such that ${\frak c} \circ g={\frak b}$.
Since $(X \times \widehat{X})^{s_0+s_1}$ is an abelian variety and
$M_H(w)$ is smooth, ${\frak c}$ extends to 
a morphism $M_H(v) \to (X \times \widehat{X})^{s_0+s_1}$.
We also denote it by ${\frak c}$.

On the other hand, ${\frak a}:M_H(v) \to Y \times \widehat{Y}$
is the Albanese map of $M_H(v)$.
Hence there is a morphism $a:Y \times \widehat{Y} \to
(X \times \widehat{X})^{s_0+s_1}$ with 
$a \circ {\frak a}={\frak c}$ and
 we have the following commutative diagram:
\begin{equation*}
\begin{CD}
M_H(v) @<{g}<< P^s\\
@V{\frak a}VV @VV{{\frak b}}V\\
Y \times \widehat{Y} @>{a}>> (X \times \widehat{X})^{s_0+s_1} 
\end{CD}
\end{equation*}
Hence $\dim \im {\frak b} \leq 4$, which implies that
$H^j(\Phi_{Y \to X}^{{\bf E}}({\cal E})[k])$, $j=0,1$ are families
of semi-homogeneous sheaves.
We set $v_i:=n_i' v_i'$, where $v_i'$ are primitive.
Then $\langle v^2 \rangle=
\langle (v_0-v_1)^2 \rangle=-2n_0' n_1' \langle v_0',v_1' \rangle$.
Hence $\langle v_0',v_1' \rangle<0$.
For a point $q \in Q(w)_0$, $V_{\bullet}:V_0 \to V_1$
denotes $({\cal V}_{\bullet})_q$.
We set $A_i:=H^i(V_{\bullet})$.
Then $A_i$ are semi-homogeneous sheaves with $v(A_i)=v_i$.
Since $\Hom_{{\bf D}(X)}(V_{\bullet},V_{\bullet}) \cong {\Bbb C}$,
$V_{\bullet}$ is not quasi-isomorphic to 
$A_0 \oplus A_1[1]$.
Hence $\Hom_{{\bf D}(X)}(A_1[-1],A_0[1]) \ne 0$.
Then 
$\Ext^2(A_1,A_0)=\Hom_{{\bf D}(X)}(A_1[-1],A_0[1]) \ne 0$ and 
$\Ext^1(A_1,A_0)=\Hom(A_1,A_0)=0$ (see \eqref{eq:aut} and
Remark \ref{rem:-1} in Appendix).
By Proposition \ref{prop:semi-homog},
$\Ext^i(({\cal A}_1)_{q_1},({\cal A}_0)_{q_0})=0$,
$i \ne 0$ for all $(q_0,q_1) \in Q_0' \times Q_1'$ and
 $\Ext^2_{p_{Q_0' \times Q_1'}}({\cal A}_1,{\cal A}_0)$ is
a locally free sheaf on $Q_0' \times Q_1'$
and all 2-extensions are parametrized by
the associated vector bundle $\overline{P} \to Q_0' \times Q_1'$, 
where we also denote the pull-backs of ${\cal A}_i$
to $Q_0' \times Q_1' \times X$ by the same ${\cal A}_i$.
$\overline{P}$ is a quotient bundle of $P$.
We denote the image of $P^s$ to $\overline{P}$ by $\overline{P}^s$.
Then we have a morphism 
$\overline{g}:\overline{P}^s \to M_H(v)$
such that $g$ is the composite 
$P^s \to \overline{P}^s \overset{\overline{g}}{\to} M_H(v)$.
Since ${\cal A}_i$ are $GL(W_i)$-equivariant, 
$G:=(GL(W_0) \times GL(W_1))/{\Bbb C}^{\times}$ acts on 
$\overline{P}$.
By Lemma \ref{lem:action} in Appendix, 
$G$ acts freely on $\overline{P}^s$
and the fiber of $\overline{g}$ is the $G$-orbit.
By Corollary \ref{cor:dim},
$\dim Q_i'-\dim GL(W_i)=\dim {\cal M}_H(n_i' v_i')^{ss}=n_i'$, and hence
\begin{equation}
\begin{split}
\dim \overline{g}(\overline{P}^s)=& \dim \Ext^2(A_1,A_0)+n_0'+n_1'+1\\
=&-n_0' n_1'\langle v_0',v_1' \rangle+n_0'+n_1'+1.
\end{split}
\end{equation}
Then we get
\begin{equation}
\begin{split}
\dim M_H(v)-\dim \overline{g}(\overline{P}^s)=&
-2n_0' n_1' \langle v_0',v_1' \rangle+2
-(-n_0' n_1'\langle v_0',v_1' \rangle+n_0'+n_1'+1)\\
= &n_0' n_1'(-\langle v_0',v_1' \rangle-1)+ (n_0'-1)(n_1'-1),
\end{split}
\end{equation}
which implies that 
$\langle v_0',v_1' \rangle=-1$ and $ (n_0'-1)(n_1'-1)=0$.
The last claim follows from \cite[Cor. 4.15]{Y:7}.
\end{proof}

\begin{rem}\label{rem:component}
We note that $g$ extends to a morphism from
an open subset of $P$.
Hence even if we do not know
$\dim \Alb(M_{H'}(w))$,
the closure of $Q_i'$ should be a union of irreducible components
of $Q_i$.
\end{rem}

\begin{rem}
In the proof of Lemma \ref{lem:induction} below,
we shall see that $M_H(v)$ is birationally equivalent
to $\widehat{Z} \times \Hilb_Z^{\langle v^2 \rangle/2}$
for an abelian surface $Z$.
\end{rem}

\subsection{Reduction to the case 
where $V_{\bullet}$ is a sheaf up to shift.}\label{subsect:reduction}
If $\rk A_0=\rk A_1=0$, then $c_1(A_0)$ and $c_1(A_1)$ are effective,
and hence
$\langle v(A_0),v(A_1) \rangle =(c_1(A_0),c_1(A_1)) \geq 0$.
This is a contradiction.
Since $\Hom(A_0,A_1)=\Ext^2(A_1,A_0)^{\vee} \ne 0$,
we see that $A_0$ is locally free.  
We first show the following:

\begin{lem}\label{lem:induction}
Keep notations as above.
There is a Fourier-Mukai functor
${\cal F}:{\bf D}(X) \to {\bf D}(X)_{op}$ such that
${\cal F}(v)=v$ and one of the following three conditions holds 
\begin{enumerate}
\item
[(1)]
 $\rk(H^0({\cal F}(V_{\bullet})))+\rk
H^1({\cal F}(V_{\bullet})))<\rk H^0(V_{\bullet})+\rk H^1(V_{\bullet})$,
or 
\item
[(2)]
 $\deg H^1({\cal F}(V_{\bullet}))<\deg H^1(V_{\bullet})$ if
$\rk H^1(V_{\bullet})=0$, or
\item
[(3)]
 $H^1({\cal F}(V_{\bullet}))=0$ if $H^1(V_{\bullet})$ is of 0-dimensional.
\end{enumerate}
\end{lem}
\begin{proof}
(i) 
We first assume that $n_1'=1$.
Since $\langle v_0',v_1' \rangle=-1$, $X_0:=M_H(v_0')$ is a fine moduli 
space.
Let ${\bf F}$ be the universal family of stable 
semi-homogeneous sheaves on $X_0 \times X$.
Applying $\Phi_{X \to X_0}^{{\bf F}^{\vee}}$ to the exact
triangle
\begin{equation}
A_0 \to V_{\bullet} \to A_1[-1] \to A_0[1],
\end{equation}
we get an exact triangle
\begin{equation}\label{eq:triangle1}
\Phi_{X \to X_0}^{{\bf F}^{\vee}}(A_0) \to 
\Phi_{X \to X_0}^{{\bf F}^{\vee}}(V_{\bullet}) \to 
\Phi_{X \to X_0}^{{\bf F}^{\vee}}(A_1)[-1] \to 
\Phi_{X \to X_0}^{{\bf F}^{\vee}}(A_0)[1].
\end{equation}
By Proposition \ref{prop:semi-homog},
$L:=\Phi_{X \to X_0}^{{\bf F}^{\vee}}({A}_1)$ is a line bundle
on $X_0$.
We note that $G:=\Phi_{X \to X_0}^{{\bf F}^{\vee}}({A}_0)[2]$ 
is a 0-dimensional sheaf
of length $n_0'$ on $X_0$. Hence \eqref{eq:triangle1}
is 
\begin{equation}\label{eq:triangle2}
G[-1] \to 
\Phi_{X \to X_0}^{{\bf F}^{\vee}}(V_{\bullet})[1] \to L 
\overset{f}{\to} G.
\end{equation}
We can choose a general point $q \in Q(w)_0$ 
such that $f$ is surjective. Then
$G \cong {\cal O}_Z \otimes L$
for a 0-dimensional subscheme $Z$ of $n_0'$-points
and we get an exact
sequence
\begin{equation}
0 \to H^1(\Phi_{X \to X_0}^{{\bf F}^{\vee}}(V_{\bullet})) \to
L \overset{f}{\to}
{\cal O}_Z \otimes L \to 0.
\end{equation}
Thus $\Phi_{X \to X_0}^{{\bf F}^{\vee}}(V_{\bullet})=I_Z \otimes L[-1]$.
By taking the dual, we have an exact triangle
\begin{equation}\label{eq:dual}
{\cal O}_Z^{\vee} \to L^{\vee} \to (I_Z \otimes L)^{\vee} \to
{\cal O}_Z^{\vee}[1].
\end{equation}
We note that ${\cal O}_Z^{\vee}=
{\cal E}xt^2_{{\cal O}_{X_0}}({\cal O}_Z,{\cal O}_{X_0})[-2]
\cong {\cal O}_Z[-2]$, if $Z$ consists of disjoint $n_1'$-points.  
We fix a line bundle $L_0$ on $X_0$ with
$c_1(L_0)=c_1(L)$.
For \eqref{eq:dual}, 
by taking a tensor product
$\otimes L_0^{\otimes 2}$ and
applying $\Phi_{X_0 \to X}^{{\bf F}}$,
we get an exact triangle
\begin{equation}
\Phi_{X_0 \to X}^{{\bf F}}(I_Z^{\vee}
 \otimes (L^{\vee} \otimes L_0^{\otimes 2}))[1]
 \to A_0 \overset{e'}{\to} B_1 \to
 \Phi_{X_0 \to X}^{{\bf F}}(I_Z^{\vee}
 \otimes (L^{\vee} \otimes L_0^{\otimes 2}))[2],
\end{equation} 
where $A_0 \cong
H^2(\Phi_{X_0 \to X}^{{\bf F}}({\cal O}_Z^{\vee} \otimes 
(L^{\vee} \otimes L_0^{\otimes 2})))$ and $B_1:= 
H^2(\Phi_{X_0 \to X}^{{\bf F}}(L^{\vee} \otimes L_0^{\otimes 2}))$
is a stable semi-homogeneous sheaf with $v(B_1)=v(A_1)$.
We set $A_0':=\ker e'$ and
$A_1':=\coker e'$. 
By shrinking $Q(w)_0$,
we may assume that $A_i'$ form flat families over $Q(w)_0$. 
Since $A_0 \to B_1$ is not trivial, we get the assertions.

(ii)
We next assume that $n_0'=1$. In this case, we consider the Fourier-Mukai 
transform $\Phi_{X \to X_1}^{{\bf F}^{\vee}}:{\bf D}(X) \to {\bf D}(X_1)$,
where $X_1:=M_H(v_1')$ and ${\bf F}$ is the universal family on
$X_1 \times X$.
Then we have an exact triangle
\begin{equation}\label{eq:n0-1}
L \to
\Phi_{X \to X_1}^{{\bf F}^{\vee}}(V_{\bullet})[2] \to 
\Phi_{X \to X_1}^{{\bf F}^{\vee}}(A_1[1]) \to
L[1]
\end{equation}
where $L:=\Phi_{X \to X_1}^{{\bf F}^{\vee}}(A_0)[2]$ is a
line bundle on $X_0$.
For a general $q \in Q(w)_0$, we may assume that
$\Phi_{X \to X_1}^{{\bf F}^{\vee}}(A_1)={\cal O}_{Z}[2]$ for 
a subscheme of distinct $n_1'$-points $Z$ on $X_1$.
Then $({\cal O}_Z[2])^{\vee} \cong {\cal O}_Z$. Hence
by taking the dual of \eqref{eq:n0-1}, we get an exact triangle
\begin{equation}
L^{\vee} \to {\cal O}_Z \to 
\Phi_{X \to X_1}^{{\bf F}^{\vee}}(V_{\bullet})^{\vee}[-1]
\to L^{\vee}[1].
\end{equation}
We fix a line bundle $L_1$ on $X_1$ with
$c_1(L_1)=c_1(L)$.
Since $B_0:=
\Phi_{X \to X_1}^{{\bf F}}(L^{\vee} \otimes L_1^{\otimes 2})$ is a
stable semi-homogeneous vector bundle with the Mukai vector
$v_1$ and
$\Phi_{X \to X_1}^{{\bf F}}({\cal O}_Z)=A_1$, we have an exact triangle
\begin{equation}
B_0 \to A_1 \to 
\Phi_{X \to X_1}^{{\bf F}}(\Phi_{X \to X_1}^{{\bf F}^{\vee}}(V_{\bullet})^{\vee} \otimes L_1^{\otimes 2})[-1]
\to B_0 [1],
\end{equation}
which implies that the assertions holds.
\end{proof}

Applying Lemma \ref{lem:induction} successively, 
we get a Fourier-Mukai functor 
${\cal F}:{\bf D}(X) \to {\bf D}(X)$ or
${\cal F}:{\bf D}(X) \to {\bf D}(X)_{op}$
such that ${\cal F}(V_{\bullet})$ is a sheaf with 
$v({\cal F}(V_{\bullet}))=v$.

\subsection{Proof of Theorem \ref{thm:main}
(the case where $\Phi_{Y \to X}^{{\bf E}}(E)$ is a sheaf)}\label{subsect:sheaf}
Replacing $V_{\bullet}$ by ${\cal F}(V_{\bullet})$,
we may assume that
$\WIT_k$ holds for a general ${\cal E}_q$.
Assume that $V:=H^k(\Phi_{Y \to X}^{{\bf E}}({\cal E}_q))$ is not stable.
By \cite[Thm. 4.14]{Y:7}, $V$ fits in an exact sequence

\begin{equation}\label{eq:HN}
0 \to A_0 \to V \to A_1 \to 0,
\end{equation}
where
$A_i$ are semi-homogeneous sheaves with the Mukai vector $n_i' v_i'$,
$\langle v_0',v_1' \rangle=1$ and $(n_0'-1)(n_1'-1)=0$.
We may assume that $A_i$ are direct sum of distinct 
stable sheaves $A_{ij} \in M_H(v_i')$, $j=1,2,\dots,n_i'$.
By using the following lemma, we shall replace the extension
\eqref{eq:HN} by an extension in another direction.

\begin{lem}\label{lem:dual}
Let $V$ fits in an exact sequence
\begin{equation}
0 \to A_0 \to V \to A_1 \to 0,
\end{equation}
with 
$A_i=\oplus_j A_{ij}$, $A_{ij} \in M_H(v_i')$, $j=1,2,\dots,n_i'$
and $\langle v_0',v_1' \rangle=1$.
Then there is a Fourier-Mukai transform
${\cal F}:{\bf D}(X) \to {\bf D}(X)_{op}$ such that
${\cal F}(V)$ fits in an exact sequence 
\begin{equation}
0 \to B_1 \to {\cal F}(V) \to B_0 \to 0,
\end{equation}
where $B_i=\oplus_j B_{ij}$, $B_{ij} \in M_H(v_i')$, $j=1,2,\dots,n_i'$.
\end{lem}

\begin{proof}
By the symmetry of the condition, we may assume that $n_1'=1$. 
We set $X_0:=M_H(v_0')$ and
${\bf F}$ the universal family on $X_0 \times X$.
Since $\chi(A_1,A_0)<0$,
$\IT_1$ holds for $A_1$ and
$L:=H^1(\Phi_{X \to X_0}^{{\bf F}^{\vee}}(A_1))$ is a line bundle on
$X_0$.
We fix a line bundle $L_0$ with $c_1(L_0)=c_1(L)$.
Then we see that
$V':=\Phi_{X_0 \to X}^{{\bf F}}
(\Phi_{X \to X_0}^{{\bf F}^{\vee}}(V)^{\vee} \otimes L_1^{\otimes 2})$
is a sheaf and fits in an exact sequence
\begin{equation}
0 \to B_1 \to V' \to A_0 \to 0,
\end{equation}
where $B_1:=\Phi_{X_0 \to X}^{{\bf F}}
(L^{\vee} \otimes L_0^{\otimes 2})[1] \in M_H(v_1)$.
We set $B_0:=A_0$. Then the claim holds.
\end{proof}

We shall show that the instability is improved, under the operation
${\cal F}$ in Lemma \ref{lem:dual}.
We only treat the case where $\rk V>0$.
The other cases are similar.
For the exact sequence \eqref{eq:HN}, by using Lemma \ref{lem:dual},
we replace $V$ by ${\cal F}(V)$.
Since \eqref{eq:HN} is the Harder-Narasimhan filtration,
$A_1$ and hence $B_1$ is locally free.
Assume that $V':={\cal F}(V)$ is not stable for all point 
$q \in Q(w)$. Then a general $V'$ fits in an exact sequence
\begin{equation}
0 \to A_0' \to V' \to A_1' \to 0,
\end{equation}
where (1) $A_i'=\oplus_j A_{ij}', i=0,1$ are direct sum
of distinct stable semi-homogeneous sheaves $A_{ij}'$ with
$v(A_{ij}')=v(A_{ik}')$ for all $j,k$, and
(2) $A_0'$ is a torsion sheaf, or
$V'$ is torsion free and $0 \subset A_0' \subset V'$ is the
Harder-Narasimhan filtration of $V'$.

We shall divide the proof into three cases
\begin{enumerate}
\item[(a)]
$V$ is not torsion free.
\item[(b)]
$V$ is torsion free but not $\mu$-semi-stable.
\item[(c)]
$V$ is $\mu$-semi-stable, but not stable.
\end{enumerate}
(a) Assume that $V$ has a torsion. Then 
$A_0$ is the torsion submodule of $V$.
Since $V$ is simple, we see that $\deg A_0>0$.
We show that the degree of the torsion submodule of $V'$
is strictly smaller than that of $V$,
that is,   
$\deg A_0'<\deg A_0$, if $V'$ has a torsion.
Assume that $V'$ has a torsion.
Then $A_0'$ is the torsion submodule of $V'$.
Since $B_1$ is locally free, 
$\varphi:A_0' \to B_0$ is injective.
If $\deg A_0'=\deg B_0$, then $\varphi$ is surjective in codimension 1.
By using the locally freeness of $B_1$, 
we see that $V' \cong B_1 \oplus B_0$,
which is a contradiction.
Thus $\deg(A_0')<\deg B_0$.

(b) Assume that $V$ is torsion free, but not $\mu$-semi-stable.
Then $B_0$ is also locally free.
If $\mu(A_0')>\mu(V')$, then
$A_0' \to B_0$ is not zero, which implies that
$\mu(A_0') \leq \mu(B_0)$.
If $\mu(A_0')=\mu(B_0)$, then we also see that 
$A_0' \to B_0$ is injective, $n_0'=1$ and
$A_0'$ is a direct summand of $V'$.
Therefore $\mu(A_0')<\mu(B_0)=\mu(A_0)$.
We can also see that $\mu(A_1')>\mu(A_1)$.
Indeed since $\mu(A_0') \geq \mu(V')>\mu(B_1)$,
$A_1' \to B_1$ is not zero.
Then we have a non-trivial homomorphism
$A_{1j}' \to B_1$ for a $j$ and we see that
$\mu(A_1')=\mu(A_{1j}')>\mu(A_1)$.

(c) If $V$ is $\mu$-semi-stable, i.e.,
$\mu(A_0)=\mu(A_1)$, then by the same argument, we see that 
$\chi(A_0')/\rk A_0'<\chi(A_0)/\rk A_0$ and
$\chi(A_1')/\rk A_1'>\chi(A_1)/\rk A_1$.
Therefore 
by applying Lemma \ref{lem:dual} successively, we get a stable sheaf.
Thus we complete the proof of Theorem \ref{thm:main}.
\qed

\subsection{In the case where $Y$ is not fine}

In the notation in section \ref{subsect:FM},
even if $Y$ is not fine, there is a universal family as a
$p_Y^*(\alpha^{-1})$-twisted sheaf for a suitable 
Cech 2-cocycle $\alpha$ of ${\cal O}_X^{\times}$.
Then we have an equivalence 
\begin{equation}
\Phi_{Y \to X}^{{\bf E}}:{\bf D}^{\alpha}(Y) \to {\bf D}(X),
\end{equation}
where ${\bf D}^{\alpha}(Y):={\bf D}(\Coh^{\alpha}(X))$ 
is the bounded derived category of coherent $\alpha$-twisted sheaves.
Let $M_{H'}^{\alpha}(w)$ be the moduli space of stable $\alpha$-twisted 
stable sheaves $E$ with $v(E)=w$.
If $\dim \Alb M_H(w)=4$, then Theorem \ref{thm:main} also holds
for this case.
By a similar method as in \cite{Y:13},
we can show that $\dim \Alb(M_{H'}^{\alpha}(w))=4$, 
if $\langle w^2 \rangle>0$
(cf. \ref{subsect:twist} in Appendix).
Here we treat one example by another argument based on \cite{Y:12}.
In the same way as in \cite[Prop. 3.14]{Y:12},
we see that for a stable $\alpha$-twisted sheaf $E$ of rank $0$,
$\Phi_{Y \to X}^{{\bf E}}(E(nH'))$ is stable for $n \gg 0$.
In particular, we have an isomorphism
$M_{H'}^{\alpha}(w e^{nH'}) \cong M_H(v')$, where $v(E)=w$ and
$v(\Phi_{Y \to X}^{\alpha}(E(nH')))=v'$.
In particular $\Alb(M_{H'}^{\alpha}(w)) \cong X \times \widehat{X}$.
Then by the same proof as in Theorem \ref{thm:main},
we see that
$M_{H'}^{\alpha}(w)$ is birationally equivalent to $M_H(v)$.
Since the support map $M_{H'}^{\alpha}(w) \to \Hilb_Y^{c_1(w)}$
($E \mapsto \Div(E)$) is a Lagrangian fibration, 
we get the following:
\begin{prop}
Assume that there is an primitive isotropic vector $v_0$
such that $v_0$ is algebraic and $\langle v,v_0 \rangle=0$.
Then $M_H(v)$ is birationally equivalent to
a holomorphic symplectic manifold with a Lagrangian fibration.
\end{prop}

\begin{cor}
The Albanese fiber $K_H(v)$ is birationally equivalent to
an irreducible symplectic manifold with a Lagrangian fibration 
if and only if $\Pic(K_H(v))$ has an isotropic element
with respect to the Beauville form.
\end{cor}
For related results on Lagrangian fibrations on 
irreducible symplectic manifolds, see
\cite{G:1},\cite{Mark:1},\cite{S:1} and references therein.

\section{Appendix}\label{sect:appendix}

\subsection{Semi-homogeneous sheaves}

The following assertions are well-known (cf. \cite{Mu:1}, \cite{Or:1}).
\begin{prop}\label{prop:semi-homog}
Let $E$ and $F$ be semi-homogeneous sheaves.
\begin{enumerate}
\item
Assume that $E$ and $F$ are locally free sheaves.
\begin{enumerate}
\item
If $\langle v(E),v(F) \rangle>0$, then
$\Hom(E,F)=\Ext^2(E,F)=0$.
\item
If $\langle v(E),v(F) \rangle<0$, then
$\mu(E) \ne \mu(F)$, $\Ext^1(E,F)=0$ and
\begin{equation}
\begin{cases}
\Hom(E,F)=0, & \mu(E)>\mu(F)\\
\Ext^2(E,F)=0, & \mu(F)>\mu(E).
\end{cases}
\end{equation}
\end{enumerate}
\item
Assume that $E$ is locally free and $F$ is a torsion sheaf.
\begin{enumerate}
\item
If $\langle v(E),v(F) \rangle>0$, then
$\Hom(E,F)=\Ext^2(E,F)=0$.
\item
If $\langle v(E),v(F) \rangle<0$, then
$\Ext^1(E,F)=\Ext^2(E,F)=0$.
\end{enumerate}
\item
Assume that $E$ and $F$ are torsion sheaves.
Then $\langle v(E),v(F) \rangle \geq 0$.
If $\langle v(E),v(F) \rangle>0$, then 
$\Hom(E,F)=\Ext^2(E,F)=0$.
\end{enumerate}
\end{prop}
This is equivalent to the fact that
the Fourier-Mukai transform of a semi-homogeneous sheaf is a sheaf
up to shift.

\begin{proof}
We only prove (i).
Indeed the proof of (ii) and (iii) are reduced to (i)
via a suitable Fourier-Mukai transform.
We note that $E^{\vee} \otimes F$ is semi-homogeneous.
There is a filtration
$ \subset F_1 \subset F_2 \subset \cdots \subset F_s=E^{\vee} \otimes F$
such that
$E_i=F_i/F_{i=1}$, $1 \leq i \leq s$ are simple semi-homogeneous vector 
bundles with $c_1(E_i)/\rk E_i=c_1(E^{\vee} \otimes F)/(\rk E \rk F)=
c_1(F)/\rk F-c_1(E)/\rk E$.
Since $\chi(E_i)/\rk E_i=(c_1(E_i)/\rk E_i)^2 /2=\chi(E,F)/\rk E \rk F$,
it is sufficient to prove the claim for $E={\cal O}_X$ and
$F$ is a simple semi-homogeneous vector bundle.
Then there is an isogeny $\pi:Y \to X$ and a line 
bundle $L$ on $Y$ such that $\pi_*(L)=F$.
Hence $\Ext^i(E,F)=H^i(X,F)=H^i(Y,L)$.
In particular, $\chi(E,F)=\chi(L)=(c_1(L)^2)/2$.
If $\chi(L)<0$, then
$H^i(Y,L)=0$ for $i \ne 1$. Thus (a) holds.
Assume that $\chi(L)>0$.
Since $\pi_*(c_1(L))=c_1(F)$,
$(c_1(L),\pi^*(H))=(c_1(F),H)$.
If $\mu(F)=0$, then the Hodge index theorem implies that
$(c_1(L)^2) \leq 0$, which is a contradiction.
Therefore $\mu(F) \ne 0=\mu(E)$. The other claims also follow
from $(c_1(L),\pi^*(H))=(c_1(F),H)$.
\end{proof}

\subsection{2-extensions}
We collect some elementary facts on 2-extensions.
We have a natural map
\begin{equation}
\Xi:\Ext^2(A_1,A_0) \to Ob({\bf D}(X))/(\text{quasi-isom.})
\end{equation}
by sending a 2-extension class
\begin{equation}
0 \to A_0 \to V_0 \to V_1 \to A_1 \to 0
\end{equation}
to the complex $V_{\bullet}:V_0 \to V_1$. 
We want to study the fiber of $\Xi$.
We take a resolution
\begin{equation}
0 \to E_{-2} \to E_{-1} \to E_0 \to A_1 \to 0
\end{equation}
such that $H^j(X,E_i^{\vee} \otimes A_0)=0$ for
$i=0,-1$, $j>0$.
Then we also have $H^j(X,E_{-2}^{\vee} \otimes A_0)=0$ for $j>0$.
Hence $\Ext^2(A_1,A_0) \cong \Hom(E_{-2},A_0)/\im(\Hom(E_{-1},A_0))$ 
and for a representative $\varphi \in  \Hom(E_{-2},A_0)$,
$\Xi([\varphi])$ is the cone $V_{\bullet}$ defined by
\begin{equation}
\varphi:E_{\bullet}[-2] \to A_0.
\end{equation}
For two exact triangles,
\begin{equation}
A_0 \to V_{\bullet}^i \to A_1[-1] \to A_0[1],\;i=1,2,
\end{equation}
we have an exact and commutative diagram:
\begin{equation}
\begin{CD}
@. 0 @. 0 @. 0 \\
@. @VVV @VVV @VVV\\
0 @>>> \Hom_{{\bf D}(X)}(A_1[-1],A_0) @>>> 
\Hom_{{\bf D}(X)}(V_{\bullet}^1,A_0) @>>>
\Hom(A_0,A_0)\\
@. @VVV @VVV @VVV\\
0 @>>> \Hom_{{\bf D}(X)}(A_1[-1],V_{\bullet}^2) 
@>>> \Hom_{{\bf D}(X)}(V_{\bullet}^1,V_{\bullet}^2) 
@>>>\Hom_{{\bf D}(X)}(A_0,V_{\bullet}^2)\\
@. @VVV @VVV @VVV\\
0 @>>> \Hom(A_1[-1],A_1[-1]) @>>> 
\Hom_{{\bf D}(X)}(V_{\bullet}^1,A_1[-1]) @>>>
\Hom_{{\bf D}(X)}(A_0,A_1[-1])=0.
\end{CD}
\end{equation}
Hence we have an exact sequence
\begin{equation}\label{eq:aut}
0 \to \Hom_{{\bf D}(X)}(A_1[-1],A_0) \overset{i}{\to}
\Hom_{{\bf D}(X)}(V_{\bullet}^1,V_{\bullet}^2) \overset{r}{\to} 
\Hom(A_0,A_0) \oplus \Hom(A_1[-1],A_1[-1]).
\end{equation}

We take a quasi-isomorphism $(V_{\bullet}^1)' \to V_{\bullet}^1$ 
such that $\Ext^j((V_1^1)',V_i^2)=\Ext^j((V_1^1)',A_0)=0$ for
$j>0,i=0,1$ and $(V_i^1)'=0$ for $i \ne 0,1$.
Then
$\Hom_{{\bf D}(X)}(V_{\bullet}^1,V_{\bullet}^2)$ 
is the cohomology group of the complex
\begin{equation}
\Hom((V_1^1)',V_0^2) \to \Hom((V_0^1)',V_0^2) \oplus 
\Hom((V_1^1)',V_1^2) \to \Hom((V_0^1)',V_1^2).
\end{equation}
Then $\varphi \in \Hom_{{\bf D}(X)}(V_{\bullet}^1,V_{\bullet}^2)$ 
induces an exact and commutative diagram:
\begin{equation}
\begin{CD}
0 @>>> A_0 @>>> (V_0^1)' @>{\phi'}>> (V_1^1)' @>>> A_1 @>>> 0\\
@. @VVV @V{\varphi_0}VV @VV{\varphi_1}V @VVV @.\\
0 @>>> A_0 @>>> V_0^2 @>{\phi}>> V_1^2 @>>> A_1 @>>> 0
\end{CD}
\end{equation}
Conversely this diagram gives an element 
$\phi \in \Hom_{{\bf D}(X)}(V_{\bullet}^1,V_{\bullet}^2)$.
For 
\begin{equation}
\varphi \in \Hom_{{\bf D}(X)}(A_1[-1],A_0)=
\Hom(\im \phi',A_0)/\Hom((V_1^1)',A_0),
\end{equation}
$i(\varphi)$ is represented by
$(\varphi \circ \phi',0) \in 
\Hom((V_0^1)',V_0^2) \oplus \Hom((V_1^1)',V_1^2)$.

We have an action of $\Aut(A_0) \times \Aut(A_1)$ on $\Ext^2(A_1,A_0)$:
\begin{equation}
\begin{matrix}
(g_0,g_1):& \Ext^2(A_1,A_0)& \to & \Ext^2(A_1,A_0)\\
& e & \mapsto & g_0 \cup e \cup g_1^{-1}.
\end{matrix}
\end{equation}
It is easy to see that the following lemma holds.
\begin{lem}\label{lem:action}
\begin{equation}
\begin{split}
\Xi^{-1}(\Xi(e)) &=(\Aut(A_0) \times \Aut(A_1))e,\\
r(\Aut_{{\bf D}(X)}(V_{\bullet}))& =
\{(g_0,g_1)|g_0 \cup e \cup g_1^{-1}=e \}
\end{split}
\end{equation}
for $e \in \Ext^2(A_1,A_0)$ with $V_{\bullet}=\Xi(e)$.
In particular, 
$GL(W_0) \times GL(W_1)/{\Bbb C}^{\times}$ acts freely on the open subscheme
of $\overline{P}$ parametrizing 
simple complexes $V_{\bullet}$,
where $\overline{P}$ is the scheme in the proof of Theorem \ref{thm:complex}.
\end{lem}

\begin{rem}\label{rem:-1}
$\Hom(A_1,A_0) \cong 
\Hom_{{\bf D}(X)}(V_{\bullet},V_{\bullet}[-1])$.
If $V_{\bullet}$ is the Fourier-Mukai transform of a sheaf
$E$, then it is 0. 
\end{rem}

\subsubsection{Some remarks on the endomorphisms of complexes}
We shall show that 
for a complex $\widehat{V}_{\bullet}$ as in section \ref{sect:pre}, 
$\Hom_{{\bf D}(X)}(\widehat{V}_{\bullet},\widehat{V}_{\bullet})$
is represented by a morphism 
$\widehat{V}_{\bullet} \to \widehat{V}_{\bullet}$
up to homotopy.
\begin{lem}
Let $V_{\bullet}:V_0 \to V_1$ be a complex.
Let $V_{\bullet}':\cdots \to V_{-1}' \to V_0' \to V_1' \to 0$
be a complex and 
$f:V_{\bullet}' \to V_{\bullet}$ a quasi-isomorphism.
Then $\Hom_{{\bf K}(X)}(V_{\bullet},V_{\bullet}) \to
\Hom_{{\bf K}(X)}(V_{\bullet}',V_{\bullet})$ is injective,
where ${\bf K}(X)$ is the homotopy category of
complexes.
In particular, $
\Hom_{{\bf K}(X)}(V_{\bullet},V_{\bullet}) \to
\Hom_{{\bf D}(X)}(V_{\bullet},V_{\bullet})$ is injective.
\end{lem}

\begin{proof}
Let $W_{\bullet}$ be the cone of $f:V_{\bullet}' \to V_{\bullet}$.
We have an exact sequence 
$\cdots \to V_0' \to V_0 \oplus V_1' \to V_1 \to 0$.
Then
it is easy to see 
$\Hom_{{\bf K}(X)}(W_{\bullet},V_{\bullet})=0$. 
So the claim follows.
\end{proof}

\begin{NB}
\begin{cor}
If 
$\Hom_{{\bf D}(X)}(V_{\bullet},V_{\bullet}) \cong {\Bbb C}$,
Then
$
\Hom_{{\bf K}(X)}(V_{\bullet},V_{\bullet}) \to
\Hom_{{\bf D}(X)}(V_{\bullet},V_{\bullet})$ is isomorphic.
\end{cor}
\end{NB}

\begin{lem}
Let $V_{\bullet}:V_0 \overset{d}{\to} V_1$ be a complex.
We set $A_i:=H^i(V_{\bullet})$.
Assume that $\Hom(V_1,V_1) \cong \Hom(V_1,A_1)$ and
$\Ext^1(V_1,A_0)=0$.
Then
$\Hom_{{\bf D}(X)}(V_{\bullet},V_{\bullet}) \cong
\Hom_{{\bf K}(X)}(V_{\bullet},V_{\bullet})$.
\end{lem}

\begin{proof}
We take a quasi-isomorphism
$f:V_{\bullet}' \to V_{\bullet}$ such that
$H^i(V_{\bullet}')=0$, $i \ne 0,1$,
$\Ext^j(V_1',V_i)=\Ext^j(V_1',A_0)=0$
for $j>0$ and $V_1' \to V_1$ is surjective.
Then $\Hom_{{\bf K}(X)}(V_{\bullet}',V_{\bullet})
\cong \Hom_{{\bf D}(X)}(V_{\bullet},V_{\bullet})$.
We note that there is an exact and commutative diagram:
\begin{equation}
\begin{CD}
@. 0 @>>> \ker f_0 @>>> \ker f_1 @>>> 0 \\
@. @VVV @VVV @VVV @VVV\\
0 @>>> A_0 @>>> V_0' @>{d'}>> V_1' @>>> A_1 @>>> 0\\
@. @V{H^0(f)}VV @V{f_0}VV @VV{f_1}V @VV{H^1(f)}V @.\\
0 @>>> A_0 @>>> V_0 @>{d}>> V_1 @>>> A_1 @>>> 0\\
@. @VVV @VVV @VVV @VVV\\
@. 0 @. 0 @. 0 @.0 @. \\
\end{CD}
\end{equation}

It is sufficient to show the surjectivity of
$\Hom_{{\bf K}(X)}(V_{\bullet},V_{\bullet}) \to
\Hom_{{\bf K}(X)}(V_{\bullet}',V_{\bullet})$.
Let $\phi:V_{\bullet}' \to V_{\bullet}$
be a morphism.
Since $f$ is a quasi-isomorphsm, we have a morphism
$a:A_1 \to A_1$ such that
$a \circ H^1(f)=H^1(\phi)$.
By our assumption,
 there is a morphism $g:V_1 \to V_1$ with a commutative diagram
\begin{equation}
\begin{CD}
V_1 @>>> A_1\\
@V{g}VV @VV{a}V\\
V_1 @>>> A_1
\end{CD}
\end{equation}
Since $\Ext^1(V_1',A_0)=0$ and 
the image of $\phi_1-g \circ f_1$ is contained in 
$d(V_0)$, we have a morphism $\lambda:V_1' \to V_0$
such that $d \circ \lambda=\phi_1-g \circ f_1$.
Replacing $\phi_1$ by $\phi_1-d \circ \lambda$
and $\phi_0$ by $\phi_0-\lambda \circ d'$,
we may assume that
$\phi_1=g \circ f_1$.
By the above diagram, we have
$d \circ \phi_0(\ker f_0)=0$, which implies that
${\phi_0}_{|\ker f_0} \in \Hom(\ker f_0,A_0)$.
Since $\Ext^1(V_1,A_0)=0$,
there is a $\lambda':V_1' \to A_0$
such that $(\phi_0-\lambda')_{|\ker f_0}=0$.
So replacing $\phi_0$ by
$\phi_0-\lambda'$, we have a morphism
$\phi_0':V_0 \to V_0$ with $\phi_0=\phi_0' \circ f_0$.
Thus $(\phi_0',g)$ gives a desired morphism
$V_{\bullet} \to V_{\bullet}$.
\end{proof}
\begin{rem}
For a complex $V_{\bullet}$, we can find 
a quasi-isomorphism $\widehat{V}_{\bullet} \to V_{\bullet}$ 
as in section \ref{sect:pre} such that $\widehat{V}_{\bullet}$ 
satisfies the assumptions
of this lemma.
\end{rem} 

\begin{rem}
By our assumption, we see that
$\Hom(V_1,d(V_0))=0$.
Then the kernel of $r$ in \eqref{eq:aut} consists of 
$(\phi_0',g)$ such that
$g=0$ and $\phi_0'$ comes from a morphism
$d(V_0) \to A_1$.
Thus 
\begin{equation}
\ker r=\Hom(d(V_0),A_0)/\Hom(V_1,A_0)=
\Hom_{{\bf D}(X)}(A_1[-1],A_0).
\end{equation}
This is compatible with \eqref{eq:aut}.
\end{rem}

\subsection{Twisted sheaves}

Let $X=\cup_i U_i$ be an analytic open covering of
$X$ and
$\alpha=\{\alpha_{ijk}\}$ a Cech 2-cocycle of
${\cal O}_X^{\times}$ representing 
a torsion element $[\alpha] \in H^2(X,{\cal O}_X^{\times})$.
Let ${\cal M}_H^{\alpha}(v)^{ss}$ be the moduli stack of
semi-stable $\alpha$-twisted sheaves of Mukai vector $v$.

\begin{lem}\label{lem:dim}
We set $v:=(0,0,n)$.
Then $\dim {\cal M}_H^{\alpha}(v)^{ss}=n$.
\end{lem}

\begin{proof}
We fix an $\alpha$-twsited vector bundle $G$ of rank $r$ on $X$. 
Let $E$ be a 0-dimensional $\alpha$-twisted sheaf of length $n$.
Then $\Hom(G,E) \otimes G \to E$ is surjective.
We set $Q:=\Quot_{G^{\oplus rn}/X}^v$.
Then ${\cal M}_H^{\alpha}(v)^{ss}$ is the quotient stack of $Q$
by the natural action of $GL(rn)$:
${\cal M}_H^{\alpha}(v)^{ss}=[Q/GL(rn)]$.
We claim that $\dim Q=(r^2 n+1)n$.
Then $\dim {\cal M}_H^{\alpha}(v)^{ss}=\dim Q-\dim GL(rn)=
(r^2n+1)n-(rn)^2=n$ and we get our lemma.
So we shall prove the claim.
We have a natural morphism
$\phi:Q \to \overline{M}_H^{\alpha}(v)$.
Since $M_H^{\alpha}(0,0,1) \cong X$,
there is a bijective morphism $\psi:S^n X \to \overline{M}_H^{\alpha}(v)$.
In order to prove the claim, it is sufficient to
show $\dim \phi^{-1}(\psi(\sum_{i=1}^s n_i P_i))=\sum_i(r^2 n n_i-1)$,
where $P_1,...,P_s$ are distinct points of $X$.
We set $Z:=\mathrm {Spec}{\Bbb C}[[x,y]]$.
Since the punctual quot-scheme 
$\Quot_{{\cal O}_Z^{\oplus l}/Z}^{m}$ is of dimension 
$lm-1$ (cf. \cite{Y:1} or \cite[Cor. 3.7]{lecture}),
we get our claim.
\end{proof}

\begin{cor}\label{cor:dim}
Let $v_0$ be a primitive Mukai vector with 
$\langle v_0^2 \rangle=0$ and $\rk v_0>0$.
Then $\dim {\cal M}_H(nv_0)^{ss}=n$.
\end{cor}

\begin{proof}
For a sufficiently large $m$,
every semi-stable sheaf $F$ with $v(F)=nv_0$
is a quotient of $\Hom({\cal O}_X(-m),F) \otimes {\cal O}_X(-m)$:
\begin{equation}
0 \to \ker \psi \to \Hom({\cal O}_X(-m),F) \otimes {\cal O}_X(-m) 
\overset{\psi}{\to} F \to 0.
\end{equation}
We set $Y:=M_H(v_0)$.
Let ${\bf E}$ be the universal family
on $X \times Y$ as a $p_Y^*(\alpha)$-twisted sheaf,
where $\alpha$ is a suitable ${\cal O}_Y^{\times}$ 
coefficient 2-cocycle and $p_Y$ is the projection.
Since $m \gg 0$, we have an exact sequence
\begin{equation}
0 \to \Phi_{X \to Y}^{{\bf E}^{\vee}}(\ker \psi)[2] 
\to \Hom(G,E) \otimes G \to E \to 0, 
\end{equation}
where $G:=\Phi_{X \to Y}^{{\bf E}^{\vee}}({\cal O}_X(-m))[2]$
and $E:=\Phi_{X \to Y}^{{\bf E}^{\vee}}(F)[2]$.
Hence we have an isomorphism
${\cal M}_H(nv_0)^{ss} \cong [Q/GL(rn)]$,
where $r=\rk G$ and $Q$ is the scheme in Lemma \ref{lem:dim}. 
Therefore we get our claim.
\end{proof}

\subsection{Weight 1 Hodge structure}\label{subsect:twist}

Let $\alpha$ be a Cech 2-cocycle of ${\cal O}_X^{\times}$
representing a $r$-torsion element of
$H^2(X,{\cal O}_X^{\times})$.
We have a homomorphism
\begin{equation}
H^2(X,{\Bbb Z}/r{\Bbb Z}) \to H^2(X,{\cal O}_X^{\times})
\end{equation}
whose image is the set of $r$-torsion elements.
We take a representative $\xi \in H^2(X,{\Bbb Z})$
such that $[\xi \mod r] \in H^2(X,{\Bbb Z}/r{\Bbb Z})$ 
maps to $[\alpha]$.
\begin{defn}
We define a weight 1 Hodge structure on $H^{odd}(X,{\Bbb Z})$
as 
\begin{equation}
\begin{split}
H^{1,0}(H^*(X,{\Bbb Z}) \otimes {\Bbb C}):= & 
e^{\xi/r}(H^{1,0}(X)\oplus H^{2,1}(X))\\ 
H^{0,1}(H^*(X,{\Bbb Z}) \otimes {\Bbb C}):=& 
e^{\xi/r}(H^{0,1}(X)\oplus H^{1,2}(X)).
\end{split}
\end{equation}
We denote this Hodge structure by
$(H^{odd}(X,{\Bbb Z}),-\frac{\xi}{r})$.
\end{defn}
Let $v$ be a primitive Mukai vector with
$\langle v^2 \rangle>0$.
Then by a similar argument as in \cite{Y:13},
we have an isomorphism
$H^{odd}(X,{\Bbb Z}) \cong H^1(M_H^{\alpha}(v),{\Bbb Z})$
preserving the Hodge structure.
Indeed, we use the surjectivity of the period map
(the period map is a double covering).
In particular, we get $\dim \Alb(M_H^{\alpha}(v))=4$.

\vspace{1pc}

{\it Acknowledgement.}
This note is motivated by a discussion with
Tom Bridgeland on his conjecture.
I would like to thank him very much.

\end{document}